\documentclass[11pt,oneside]{article}

\usepackage{amsfonts}
\usepackage[mathscr]{euscript}
\usepackage{xcolor}

\newcommand{\lbl}[1]{\label{#1}}

\newtheorem{theo}{Theorem}
\newtheorem{prop}{Proposition}

\newtheorem{col}{Corollary}

\newcommand{\be}{\begin{equation}}
\newcommand{\ee}{\end{equation}}
\newcommand\bes{\begin{eqnarray}} \newcommand\ees{\end{eqnarray}}
\newcommand{\bess}{\begin{eqnarray*}}
\newcommand{\eess}{\end{eqnarray*}}

\newcommand\ep{\varepsilon}
\newcommand\kk{\left}
\newcommand\rr{\right}
\newcommand\dd{\displaystyle}

\newcommand\ud{\underline}
\newcommand\vp{\varphi}

\newcommand\tty{_{t\to\yy}}

\newcommand\yy{\infty}

\begin{document}
\setlength{\baselineskip}{16pt} \pagestyle{myheadings}

\begin{center}{\Large\bf Note on a two-species competition-diffusion model}\\[2mm]
{\Large\bf  with two free boundaries}\footnote{This work was
supported by NSFC Grant 11371113}\\[4mm]
 {\Large  Mingxin Wang\footnote{{\sl E-mail}: mxwang@hit.edu.cn}}\\[0.5mm]
{\small Natural Science Research Center, Harbin Institute of Technology, Harbin 150080, PR}
\end{center}

\begin{quote}
\noindent{\bf Abstract.} In \cite{GW15, Wu15}, the authors studied a two-species competition-diffusion model with two free boundaries. The existence, uniqueness and long time behavior of global solution were established. In this note we still discuss the long time behavior of global solution and provide some new results and simpler proofs.


 \end{quote}

 \section{Introduction}
 \setcounter{equation}{0} {\setlength\arraycolsep{2pt}

Recently, Guo \& Wu \cite{GW15} studied the existence and uniqueness of global solution $(u,v,s_1,s_2)$ to the following free boundary problem
  \bess
 \left\{\begin{array}{lll}
 u_t-d_1u_{xx}=r_1u(1-u-kv), &t>0,\ \ 0<x<s_1(t),\\[0.5mm]
 v_t-d_2v_{xx}=r_2v(1-v-hu),\ \ &t>0, \ \ 0<x<s_2(t),\\[0.5mm]
 u_x(t,0)=v_x(t,0)=0,\ \ \ &t\ge0, \\[0.5mm]
 s_1'(t)=-\mu_1 u_x(t,s_1(t)), \ \ s_2'(t)=-\mu_2 v_x(t,s_2(t)), \ &t\ge0,\\[0.5mm]
 u=0 \ {\rm for}\ x\ge s_1(t), \ \ v=0 \ {\rm for}\ x\ge s_2(t), &t\ge0,\\[0.5mm]
 u(0,x)=u_0(x), \ \ v(0,x)=v_0(x), \ &x\in[0, \yy), \\[0.5mm]
 s_1(0)=s_1^0>0, \ \ s_2(0)=s_2^0>0,
 \end{array}\right.
 \eess
where the parameters are positive constants, and $u_0(x),v_0(x)$ satisfy
 \bess
 &u_0\in C^2([0,s_1^0]),\  u'_0(0)=0,\  u_0(x)>0\ \ {\rm in}\, \
[0,s_1^0), \ u_0(x)=0\ \ {\rm in}\, \ [s_1^0,\yy),&\\[0.5mm]
 & v_0\in C^2([0,s_2^0]),\  v'_0(0)=0,\  v_0(x)>0\ \ {\rm in}\, \
[0,s_2^0), \ v_0(x)=0\ \ {\rm in}\, \ [s_2^0,\yy).&
 \eess
Furthermore, Guo \& Wu \cite{GW15} and  Wu \cite{Wu15} investigated the long time behavior of $(u,v,s_1,s_2)$ for the cases $0<k<1<h$ and $0<k,\,h<1$, respectively.

By use of the arguments of \cite[Theorem 2.1]{Wperiodic15} we can prove that
$ s_1'(t), s_2'(t)>0$, and
 \bess
    (u,v,s_1,s_2)\in C^{1+\frac{\alpha}2,2+\alpha}({\cal D}_\yy^{s_1})\times
 C^{1+\frac{\alpha}2,2+\alpha}({\cal D}_\yy^{s_2})\times [C^{1+\frac{1+\alpha}2}([0,\yy))]^2,
 \eess
where ${\cal D}_\yy^{s_i}=\{t\ge 0,\,0\le x\le s_i(t)\}$. Moreover, there exists a positive constant $C$ such that
  \bes\left\{\begin{array}{ll}
  \|u(t,\cdot)\|_{C^1([0,\,s_1(t)])}, \ \ \|v(t,\cdot)\|_{C^1([0,s_2(t)])}\leq C, \ \ &\forall \ t\ge
  1, \\[1.5mm]
  \|s_1',\,s_2'\|_{C^{{\alpha}/2}([n+1,n+3])}\leq C, \ \ &\forall \ n\geq 0.
  \end{array}\right.\lbl{a.1}\ees

We still study the long time behavior of $(u,v,s_1,s_2)$ and provide some new results and simpler proofs. This short paper can be considered as the supplements of papers  \cite{GW15, Wu15}.

\section{Preliminaries}

\begin{prop}{\rm (\cite[Proposition 2.1]{WZjdde})}\label{p1.1} \ Let $d, r,a$ be fixed positive constants. For any given $\varepsilon, L>0$, there exists $l_\varepsilon>\max\big\{L,\frac{\pi}2\sqrt{d/(ra)}\,\big\}$ such that, when a non-negative $C^{1,2}$ function $z$ satisfies
 $$ \left\{\begin{array}{ll}
  z_t-dz_{xx}\geq rz(a-z), \ \ &t>0, \ \ 0<x<l_\varepsilon,\\[.2mm]
 z_x(t,0)=0, \ \ z(t, l_\varepsilon)\geq 0,\ \ &t\ge 0
 \end{array}\right.$$
and $z(0,x)>0$ in $(0,l_\varepsilon)$, then $\liminf_{t\to\infty}z(t,x)\geq a-\varepsilon$ uniformly on $[0,L]$.
 \end{prop}

\begin{prop}\lbl{p1.2}{\rm(\cite{BDK, DLin})} For any given $d,a,b,\mu>0$, the problem
 \bes\left\{\begin{array}{l}
 d q''-c q'+q(a-b q)=0,\ \ 0<y<\infty,\\[1mm]
q(0)=0,\ \ q'(0)=c/\mu, \ \ q(\infty)=a/b,\\[1mm]
 c\in(0,2\sqrt{a d}); \ \ q'(y)>0,\ \ 0<y<\infty
\end{array}\right.\lbl{a.2}
  \ees
has a unique solution $(q,c)$. Denote $\gamma=(\mu,a,b,d)$ and $c=c(\gamma)$. Then $c(\gamma)$ is strictly increasing in $\mu$ and $a$, respectively, and is strictly decreasing in $b$. Moreover,
 \bes
 \lim_{\frac{a\mu}{b d}\to\yy}\frac{c(\gamma)}{\sqrt{a d}}=2, \ \ \
 \lim_{\frac{a\mu}{b d}\to 0}\frac{c(\gamma)}{\sqrt{a d}}\frac{b d}{a\mu}=\frac 1{\sqrt{3}}.
 \lbl{a.3}\ees
\end{prop}

To simplify the notations, we define
 \[c_1(\mu,a)=c(\mu,a,r_1,d_1), \ \ \ c_2(\mu,a)=c(\mu,a,r_2,d_2).\]
If $0<k<1$, in view of (\ref{a.3}), it is easy to see that
 \[\lim_{\mu_1\to\yy}c_1(\mu_1,r_1(1-k))=2\sqrt{d_1r_1(1-k)}, \ \ \ \lim_{\mu_2\to 0}c_2(\mu_2,r_2)=0.\]
By the monotonicity of $c(\gamma)$ in $\mu$, there exist $\mu_1^*,\mu_2^*>0$ such that $c_1(\mu_1,r_1(1-k))>c_2(\mu_2,r_2)$ for all $\mu_1\ge\mu_1^*$ and $0<\mu_2\le\mu_2^*$. Therefore, $[\mu_1^*,\yy)\times(0,\mu_2^*]\subset{\cal A}$, where
 \bes
 {\cal A}=\{(\mu_1,\mu_2):\,\mu_1,\mu_2>0,\, c_1(\mu_1,r_1(1-k))>c_2(\mu_2,r_2)\}.
 \lbl{b.4}\ees

Same as \cite{GW15, Wu15}, we define $s_i^\yy=\dd\lim_{t\to\yy}s_i(t)$, $i=1,2$, and
 \[s_1^*=\frac\pi 2\sqrt{\frac{d_1}{r_1}}, \ \ s_2^*=\frac\pi 2\sqrt{\frac{d_2}{r_2}},
 \ \ \tilde s_1=\frac\pi 2\sqrt{\frac{d_1}{r_1(1-k)}} \ \ {\rm if} \ k<1, \ \ \tilde s_2=\frac\pi 2\sqrt{\frac{d_2}{r_2(1-h)}} \ \ {\rm if} \ h<1.\]

In order to convenient writing, for any given constant $\tau\ge 0$ and function $f(t)$, we set
 \[D_\tau^f=\{(t,x):\, t\ge \tau, \ 0\le x\le f(t)\}.\]

\section{Main results and their proofs}

Using the estimates (\ref{a.1}) and \cite[Lemma 3.1]{Wjde15}, we have

 \begin{theo}\lbl{th1.1} If $s_1^\yy<\yy\ (s_2^\yy<\yy)$. Then
 \bess
 \lim_{t\to\infty}\max_{[0, s_1(t)]}u(t,\cdot)=0 \ \ \Big(\lim_{t\to\infty}\max_{[0, s_2(t)]}v(t,\cdot)=0\Big).\eess
\end{theo}

\begin{theo}\lbl{th1.2} If $s_2^\yy<\yy$, $s_1^\yy=\yy$ {\rm(}$s_1^\yy<\yy$, $s_2^\yy=\yy${\rm)}, then $\lim_{t\to\infty}u(t,x)=1$ {\em(}$\lim_{t\to\infty}v(t,x)=1${\rm)} uniformly in any compact subset of $[0,\yy)$.
 \end{theo}

{\bf Proof}. Let $z(t)$ be the unique solution of
 \[z'=r_1z(1-z), \ \ t>0; \ \ \ z(0)=\|u_0\|_{L^\yy}.\]
Then $z(t)\to 1$ as $t\to\yy$. The comparison principle leads to
 \bes
 \limsup_{t\to\infty}u(t,x)\le 1\ \ {\rm uniformly\ in} \ [0,\yy).
 \lbl{a.4}\ees
For any given $0<\ep,\delta\ll 1$ and $L>0$, let $l_\varepsilon$ be given by Proposition \ref{p1.1} with $d=d_1$, $r=r_1$ and $a=1-k\delta$. Since $s_1^\yy=\yy$, $\lim_{t\to\infty}\max_{[0, s_2(t)]}v(t,\cdot)=0$ (Theorem \ref{th1.1}) and $v=0$ for $x>s_2(t)$, there exists $T\gg 1$ such that $s_1(t)>l_\ep$ and $v(t,x)<\delta$ in $[T,\yy)\times[0,\yy)$. Thus, $u$ satisfies
 \bess
 \left\{\begin{array}{lll}
 u_t-d_1u_{xx}\ge r_1u(1-k\delta-u), &t\ge T,\ \ 0<x<l_\ep,\\[0.5mm]
 u_x(t,0)=0,\ \ u(t,l_\ep)\ge 0,\ &t\ge T.
 \end{array}\right.
 \eess
In view of Proposition \ref{p1.1}, we have $\liminf_{t\to\infty}u(t,x)\ge 1-k\delta-\ep$ uniformly on $[0,L]$. The arbitrariness of $L, \ep$ and $\delta$ imply that  $\liminf_{t\to\infty}u(t,x)\ge 1$ uniformly in any compact subset of  $[0,\yy)$. Remember  (\ref{a.4}), the desired result is obtained. \ \ \ \fbox{}

Utilizing the iteration methods used in the proof of \cite[Theorem 2.4]{WZjdde} we can prove the following theorem and the details will be omitted.

\begin{theo}\lbl{th1.3} Assume $s_1^\infty=s_2^\yy=\infty$. For any given $L>0$, the following hold:

{\rm(i)} if $0<h,k<1$, then $\lim_{t\to\infty}u(t,x)=\frac{1-k}{1-hk}$, $\lim_{t\to\infty}v(t,x)=\frac{1-h}{1-hk}$ uniformly in $[0,L]$;

{\rm(ii)}\, if $0<h<1\leq k$, then $\lim_{t\to\infty}u(t,x)=0$, $\lim_{t\to\infty}v(t,x)=1$ uniformly in $[0,L]$;

{\rm(iii)}\, if $0<k<1\leq h$, then $\lim_{t\to\infty}u(t,x)=1$, $\lim_{t\to\infty}v(t,x)=0$ uniformly in $[0,L]$.
 \end{theo}


\begin{theo}\lbl{th1.4} {\rm(i)} If $s_1^*<s_1^\yy<\yy$, then $s_2^\yy=\yy$. If $s_2^*<s_2^\yy<\yy$, then $s_1^\yy=\yy$. As the consequence, $s_1^\yy<\yy$ and $s_2^\yy>s_2^*$ imply $s_2^\yy=\infty$,
$s_2^\yy<\yy$ and $s_1^\yy>s_1^*$ imply $s_1^\yy=\infty$;

{\rm(ii)} If $k<1$ and $s_1^\yy>\tilde s_1$, then $s_1^\yy=\yy$. If $h<1$ and $s_2^\yy>\tilde s_2$, then $s_2^\yy=\yy$.
\end{theo}

{\bf Proof}. (i) We only prove the first conclusion. Because of $s_1^\yy>s_1^*$, there exist $0<\ep\ll 1$ and $T\gg 1$ such that
$1-k\ep>0$ and $s_1(t)>\frac\pi 2(\frac{d_1}{r_1(1-k\ep)})^{1/2}$ for all $t\ge T$. Assume on the contrary that $s_2^\yy<\yy$. Then $\lim_{t\to\infty}\max_{[0,s_2(t)]}v(t,\cdot)=0$ by Theorem \ref{th1.1}. There exists $T_1>T$ such that $v(t,x)<\ep$ in $D^{s_2}_{T_1}$. Therefore, $(u,s_1)$ satisfies
 \bess
 \left\{\begin{array}{lll}
 u_t-d_1u_{xx}\ge r_1u(1-u-k\ep), &t>T_1,\ \ 0<x<s_1(t),\\[0.5mm]
 u_x(t,0)=u(t,s_1(t))=0, \ &t\ge T_1,\\[.5mm]
 s_1'(t)=-\mu_1 u_x(t,s_1(t)), \ &t\ge T_1.
 \end{array}\right.
 \eess
Since $s_1(T_1)>\frac\pi 2(\frac{d_1}{r_1(1-k\ep)})^{1/2}$, we have $s_1^\yy=\yy$ (\cite[Theorem 3.4]{DLin}). A contradiction.

(ii) The proof is similar to that of (i) since the estimate (\ref{a.4}) holds true for $v$. Please refer to the proof of the following Theorem \ref{th1.5} for details. 
\ \ \ \fbox{}

 \begin{col}\lbl{c3.1} {\rm(i)} If one of $s_1^0\ge s_1^*$ and $s_2^0\ge s_2^*$ holds, then either $s_1^\yy=\yy$ or $s_2^\yy=\yy$;

 {\rm(ii)} If $k,h<1$ and $s_1^0\ge\tilde s_1$, $s_2^0\ge\tilde s_2$, then $s_1^\yy=s_2^\yy=\yy$.
 \end{col}

Theorem \ref{th1.4} and Corollary \ref{c3.1} are the improvements of \cite[Theorem 2]{GW15} and \cite[Theorem 1, Propositions 4 and 5]{Wu15}, and our proofs are very simple.

To facilitate writing, for $\tau\ge 0$, we introduce the following free boundary problem
\bes
 \left\{\begin{array}{lll}
 w_t-dw_{xx}=rw(a-w), &t>\tau,\ \ 0<x<g(t),\\[.2mm]
 w_x(t,0)=0,\ \ w(t,g(t))=0,\ \ &t\ge \tau,\\[.2mm]
 g'(t)=-\mu w_x(t,g(t)),&t\ge \tau,\\[.2mm]
 g(\tau)=g_0,\ \ w(\tau,x)=w_0(x),&0\le x\le g_0,
 \end{array}\right.\label{a.5}
 \ees
and set $\Lambda=(\tau,d,r,a,\mu,g_0)$.

\begin{theo}\lbl{th1.5} Let $0<k<1<h$ and $(\mu_1,\mu_2)\in{\cal A}$, where ${\cal A}$ is given by $(\ref{b.4})$. If $s_1^\yy>\tilde s_1$, then $s_1^\yy=\yy$ and $s_2^\yy<\yy$.
\end{theo}

{\bf Proof}. Note that $k<1$, $s_1^\yy>\tilde s_1$ and $(\mu_1,\mu_2)\in{\cal A}$, there exist $0<\ep_0\ll 1$ and $t_0\gg 1$ such that $s_1(t)>\frac\pi 2(\frac{d_1}{r_1a_\ep})^{1/2}$ and
 \[ k(1+\ep)<1,\ \ c_1(\mu_1,r_1a_\ep)>c_2(\mu_2,r_2)\]
for all $t\ge t_0$ and $0<\ep\le\ep_0$, where $a_\ep=1-k(1+\ep)$. Since the estimate (\ref{a.4}) holds true for $v$, for each fixed $0<\ep\le\ep_0$, there exists $ t_1>t_0$ such that $v(t,x)<1+\ep$ in $[t_1,\yy)\times[0,\yy)$.

Let $(w_1,g_1)$ be the unique solution of (\ref{a.5}) with $\Lambda=(t_1,d_1,r_1,a_\ep,\mu_1,s_1(t_1))$ and $w_0(x)=u(t_1,x)$.
Then $s_1(t)\ge g_1(t)$, $u(t,x)\ge w_1(t,x)$ in $D_{t_1}^{g_1}$, and $g_1(\yy)=\yy$ since $g_1(t_1)>\frac\pi 2(\frac{d_1}{r_1a_\ep})^{1/2}$. Consequently, $s_1^\yy=\yy$. Take advantage of \cite[Theorem 3.1]{ZhaoW15}, it is deduced that, as $t\to\yy$,
 \bes
 g_1(t)-\tilde c t\to\tilde\rho\in\mathbb{R}, \ \ \|w_1(t,x)- \tilde q(\tilde c t+\tilde\rho-x)\|_{L^\infty([0,g_1(t)])}\to 0,\lbl{a.8}\ees
where $(\tilde q, \tilde c)$ is the unique solution of (\ref{a.2}) with $\gamma=(\mu_1,r_1a_\ep,r_1,d_1)$, i.e., $\tilde c=c_1(\mu_1,r_1a_\ep)$.

Assume on the contrary that $s_2^\yy=\yy$. We first prove
 \bes
  \lim_{t\to\yy}\max_{[0,\yy)}v(t,\cdot)=0.\lbl{a.14}
 \ees
Let $(w,g)$ be the unique solution of (\ref{a.5}) with $\Lambda=(0,d_2,r_2,1,\mu_2,s_2^0)$ and $w_0(x)=v_0(x)$.
 Then $w(t,x)\ge v(t,x)$, $g(t)\ge s_2(t)$ in $D_0^{s_2}$. Hence, $g(\yy)=\yy$. It follows from \cite[Theorem 3.1]{ZhaoW15} that
 \bes
 \lim_{t\to\infty}(g(t)-c t)=\rho\in\mathbb{R}, \ \ \ {\rm with} \ c=c_2(\mu_2,r_2).\lbl{a.7}\ees
Thanks to $\tilde c>c$, $s_2(t)\le g(t)$ and (\ref{a.7}), it deduces that, as $t\to\yy$, $g_1(t)-g(t)\to\yy$ and
 \[\min_{x\in[0, s_2(t)]}(\tilde c t+\tilde\rho-x)\ge \tilde c t+\tilde\rho-g(t)=(\tilde c-c)t+\tilde\rho-\rho+o(1)\to\yy.\]
Based on $\tilde q(y)\nearrow a_\ep$ as $y\to\yy$, we have  
$\lim_{t\to\yy}\min_{x\in[0, s_2(t)]}\tilde q(\tilde c t+\tilde\rho-x)=a_\ep$. It then follows, upon using (\ref{a.8}), that $\lim_{t\to\yy}\min_{[0, s_2(t)]}w_1(t,\cdot)=a_\ep$. This implies $\liminf_{t\to\yy}\min_{[0, s_2(t)]}u(t,\cdot)\ge a_\ep$ since  $g_1(t)>g(t)\ge s_2(t)$ and $u(t,x)\ge w_1(t,x)$ in $D_T^{g_1}$ when $T\gg 1$. The arbitrariness of $\ep$ gives
  \bes
 \liminf_{t\to\yy}\min_{[0, s_2(t)]}u(t,\cdot)\ge 1-k.
  \lbl{a.9}\ees

Take $t_2>t_1$ such that $u(t,x)\ge 1-k-\ep$ in $D_{t_2}^{s_2}$. Let $\tilde a_\ep=1-h(1-k-\ep)$ and $(w_2,g_2)$ be the unique solution of (\ref{a.5}) with $\Lambda=(t_2,d_2,r_2,\tilde a_\ep,\mu_2,s_2(t_2))$ and $w_0(x)=v(t_2,x)$.
Then $s_2(t)\le g_2(t)$, $v(t,x)\le w_2(t,x)$ in $D_{t_2}^{s_2}$. So, $g_2(\yy)=\yy$. Similarly to the above,
\bess
 g_2(t)-c_2 t\to\rho_2\in\mathbb{R}, \ \ \|w_2(t,x)-q_2(c_2 t+\rho_2-x)\|_{L^\infty([0,\,g_2(t)])}\to 0\eess
as $t\to\infty$, where $(q_2, c_2)$ is the unique solution of (\ref{a.2}) with $\gamma=(\mu_2,r_2\tilde a_\ep,r_2, d_2)$. Then $q_2(y)<\tilde a_\ep$, and
$w_2(t,x)<\tilde a_\ep+\ep$ in $D_{t_3}^{g_2}$ for some $ t_3> t_2$. Note $s_2(t)\le g_2(t)$ for $t\ge t_3$ and $D_{t_3}^{s_2}\subset D_{t_3}^{g_2}$, the following holds:
   \bess
   v(t,x)\le w_2(t,x)<\tilde a_\ep+\ep=1-h(1-k)+(1+h)\ep \ \ {\rm in} \ D_{t_3}^{s_2}.\eess
In view of $v(t,x)=0$ for $x\ge s_2(t)$, and the arbitrariness of $\ep$, it follows that
   \bess
  \limsup_{t\to\yy}\max_{[0,\yy)}v(t,\cdot)\le 1-h(1-k):=\bar v_2 .\eess

When $h(1-k)\ge 1$, we have $\bar v_2\le 0$, and so (\ref{a.14}) holds since $v(t,x)\ge 0$.

Here we deal with the case $h(1-k)<1$. There exists $t_4\gg1$ such that $v(t,x)\le \bar v_2+\ep:=\bar v_2^\ep<1$ in $[t_4,\yy)\times[0,\yy)$. Obviously, $c_1(\mu_1,r_1(1-k\bar v_2^\ep))>c_2(\mu_2,r_2)$. Similarly to the above, we can get
 \bess
 &\dd\liminf\tty\min_{[0,s_2(t)]}u(t,\cdot)\ge 1-k[1-h(1-k)]:=\ud{u}_2,&\\[.5mm]
 &\dd \limsup_{t\to\yy}\max_{[0,\yy)}v(t,\cdot)\le 1-h\ud{u}_2:=\bar v_3.&\eess
If $h\ud{u}_2\ge 1$, then $\bar v_3=0$ and (\ref{a.14}) holds. If $h\ud{u}_2<1$, repeating the above procedure in the way of the proof of \cite[Lemma 2.2]{WZjdde}, we can get  (\ref{a.14}) eventually.

For any given $0<\delta\ll 1$, there exists $ t_5\gg 1$ such that
 $v(t,x)<\delta$ in $[t_5,\yy)\times[0,\yy)$. 
Obviously, $c_1(\mu_1,r_1(1-k\delta))>c_2(\mu_2,r_2)$. Replacing $1+\ep$ by $\delta$, similarly to the above we can prove $u(t,x)>1-\delta$ in $D_{t_6}^{s_2}$ for some $t_6>t_5$. Therefore, $1-v-hu<1-h(1-\delta)-v<0$ in $D_{t_6}^{s_2}$ because of $0<\delta\ll 1$ and $h>1$. According to \cite[Lemma 3.2]{HW}, $s_2^\yy<\yy$ is followed. This is a contradiction and the proof is finished. \ \ \ \fbox{}

Theorem \ref{th1.5} is exactly \cite[Theorem 3]{GW15}, and our proof is simpler.

Similarly to the proof of \cite[ Lemma 2.1]{WZ12}, it can be shown that
 \bess
 0<s_1(t)\le K\mu_1t+s_1^0, \ \ \forall \ t>0,
 \eess
where
 \[K=2\max\big\{\max\{1,\,\|u_0\|_{\yy}\}\sqrt{r_1/(2d_1)},\ -\min_{[0,\,s_1^0]}u_0'(x)\big\}.\]

\begin{theo}\lbl{th1.6} Let $d_i, r_i, k, h$ and $\mu_2$ be fixed. Then there exists $0<\bar\mu_1<\sqrt{2d_2r_2}/K$ such that, when
  \[0<\mu_1<\bar\mu_1, \ \ s_2^0-s_1^0> \pi\frac{2d_2}{\sqrt{2d_2r_2-K^2\mu_1^2}}:=L(\mu_1),\]
we have $s_2(t)\ge K\mu_1 t+s_1^0+L(\mu_1)\to\yy$ as $t\to\yy$. Moreover, if $k<1$ and $s_1^0\ge\tilde s_1$, we also have $s_1^\yy=\yy$ for all $\mu_1>0$.
 \end{theo}

{\bf Proof}. Denote $\sigma=K\mu_1$. For the given $\sigma\in(0,\sqrt{2d_2r_2})$, and these $t$ satisfying $s_2(t)>\sigma t+s_1^0$, we define
 \[y=x-\sigma t-s_1^0, \ \ w(t,y)=v(t,x), \ \ \eta(t)=s_2(t)-\sigma t-s_1^0.\]
Note that $y\ge 0$ implies $x\ge s_1(t)$ and $u(t,x)=0$ for $x\ge s_1(t)$, we have
 \bess\kk\{\begin{array}{ll}
 w_t-d_2w_{yy}-\sigma w_y=r_2w(1-w), \ \ &t>0, \ 0<y<\eta(t),\\[.5mm]
 w(t,0)=v(t,\sigma t+s_1^0), \ \ w(t,\eta(t))=0, \ \ &t\ge 0,\\[.5mm]
w(0,y)=v_0(0,y+s_1^0), \ \ \ &0\le y\le s_2^0-s_1^0,
 \end{array}\rr.\eess
and $w(t,y)>0$ for $t\ge 0$ and $0<y<\eta(t)$. Let $\lambda$ be the principal eigenvalue of
 \bes\left\{\begin{array}{ll}
 -d_2\phi''-\sigma\phi'-r_2\phi=\lambda\phi, \ \ 0<x<\ell,\\[1mm]
 \phi(0)=0=\phi(\ell).
 \end{array}\right.\lbl{a.10}\ees
The following relation between $\lambda$ and $\ell$ holds:
 \[\frac{\pi}\ell=\frac{\sqrt{4d_2(r_2+\lambda)-\sigma^2}}{2d_2}.\]
Take $\lambda=-r_2/2$ and define
 \[\ell_\sigma=\pi\frac{2d_2}{\sqrt{2d_2r_2-\sigma^2}}, \ \ \ \phi(y)={\rm e}^{-\frac\sigma{2d_2}y}\sin\frac{\pi}{\ell_\sigma}y.\]
Then $(\ell_\sigma,\phi)$ satisfies (\ref{a.10}) with $\lambda=-r_2/2$ and $\ell=\ell_\sigma$.
Assume $s_2^0-s_1^0>\ell_\sigma$. Set
 \[\delta_\sigma=\min\kk\{\inf_{(0,\,\ell_\sigma)}\frac{w(0,y)}{\phi(y)}, \ \frac 12\inf_{(0,\,\ell_\sigma)}\frac 1{\phi(y)}\rr\}, \ \ \psi(y)=\delta_\sigma\phi(y).\]
Then $0<\delta_\sigma<\yy$. It is easy to see that $\psi(y)\le w(0,y)$ in $[0,\ell_\sigma]$ and satisfies
   \bess\left\{\begin{array}{ll}
 -d_2\psi''-\sigma\psi'\le r_2\psi(1-\psi), \ \ 0<x<\ell_\sigma,\\[1mm]
 \psi(0)=0=\psi(\ell_\sigma).
 \end{array}\right.\eess
Take a maximal $\bar\sigma\in(0,\sqrt{2d_2r_2})$ so that
  \bes\sigma<\mu_2\delta_\sigma\frac{\pi}{\ell_\sigma}
  \exp\kk(-\frac{\sigma\ell_\sigma}{2d_2}\rr), \ \ \forall \ \sigma\in(0,\bar\sigma).\lbl{a.11}\ees

For any given $\sigma\in(0,\bar\sigma)$, we claim that $\eta(t)>\ell_\sigma$ for all $t\ge 0$, which implies $s_2(t)\ge \sigma t+s_1^0+\ell_\sigma\to\yy$. In fact, note $\eta(0)=s_2^0-s_1^0>\ell_\sigma$, if our claim is not true, then we can find a $t_0>0$ such that $\eta(t)>\ell_\sigma$ for all $0\le t<t_0$ and $\eta(t_0)=\ell_\sigma$. Therefore, $\eta'(t_0)\le 0$, i.e, $s_2'(t_0)\le \sigma$. On the other hand, by the comparison principle, we have  $w(t,y)\ge\psi(y)$ in $[0, t_0]\times[0,\ell_\sigma]$. Particularly, $w(t_0,y)\ge\psi(y)$ in $[0,\ell_\sigma]$. Due to $w(t_0,\ell_\sigma)=0=\psi(\ell_\sigma)$, one has
 $$w_y(t_0,\eta(t_0))\le\psi'(\ell_\sigma)=-\delta_\sigma\frac{\pi}{\ell_\sigma}
 \exp\kk(-\frac{\sigma\ell_\sigma}{2d_2}\rr).$$
It follows, upon using $v_x(t_0,s_2(t_0))=w_y(t_0,\eta(t_0))$, that
 \[\sigma\ge s_2'(t_0)=-\mu_2w_y(t_0,\eta(t_0))\ge\mu_2\delta_\sigma\frac{\pi}{\ell_\sigma}
 \exp\kk(-\frac{\sigma\ell_\sigma}{2d_2}\rr).\]
It is in contradiction with (\ref{a.11}).

Take $\bar\mu_1=\bar\sigma/K$. Then $0<\mu_1<\bar\mu_1$ is equivalent to $0<\sigma<\bar\sigma$.

At last, if $k<1$ and $s_1^0\ge\tilde s_1$, then $s_1^\yy=\yy$ for any $\mu_1>0$ by Theorem \ref{th1.4}(ii). The proof is complete. \ \ \fbox{}

Theorem \ref{th1.6} can be regarded as an improvement of \cite[Theorem 5]{GW15}, here we need neither the assumption $v_0'(x)\le 0$ in $[s_1^0,s_2^0]$, nor the condition that $d_2$ is suitably large. Moreover, our proof of Theorem \ref{th1.6} is simpler.

From the proof of Theorem \ref{th1.6} it can be seen that if we take $\bar\sigma\in(0,\sqrt{d_2r_2})$ such that (\ref{a.11}) holds, then $s_2^\yy=\yy$ is still true provided that $0<\mu_1<\bar\mu_1$ and $s_2^0-s_1^0\ge \pi\frac{2d_2}{\sqrt{d_2r_2}}$.

\vskip 4pt Theorem \ref{th1.5} demonstrates that when the superior competitor spreads quickly and the inferior competitor spreads slowly, the inferior competitor will vanish eventually and the superior competitor will spread successfully and occupy the whole space.
Take $0<k<h<1$ in Theorem \ref{th1.6}, the conclusion indicates that if the superior competitor spreads too slow to catch up with the inferior competitor, it may leave enough space for the inferior competitor to survive.

In the following we will discuss the more accurate limits of $(u,v)$ as $t\to\yy$ when $s_1^\yy=s_2^\yy=\yy$. By the comparison principle and \cite[Theorem 4.2]{DLin}, it can be deduced that
 \bes\kk\{\begin{array}{ll}
 \dd\liminf_{t\to\yy}(s_1(t)/t)\ge c_1(\mu_1,r_1(1-k)):=\ud{c}_1\ \ \ {\rm if} \ \ k<1,\\[1mm]
 \dd\liminf_{t\to\yy}(s_2(t)/t)\ge c_2(\mu_2,r_2(1-h)):=\ud{c}_2\ \ \ {\rm if} \ \ h<1.\end{array}\rr.
 \lbl{a.12}\ees

The following two theorems are the improvements of Theorem \ref{th1.3}.

\begin{theo}\lbl{th1.7} Let $d_i,r_i,\mu_i, k, h$ be fixed and $0<k,h<1$. If $s_1^\infty=s_2^\yy=\infty$, then
 for each $0<c_0<\min\{\ud{c}_1,\,\ud{c}_2\}$,
  \bess
 \lim_{t\to\yy}\max_{[0,\,c_0t]}\kk|u(t,\cdot)-(1-k)/(1-hk)\rr|=0, \ \ \
 \lim_{t\to\yy}\max_{[0,\,c_0t]}\kk|v(t,\cdot)-(1-h)/(1-hk)\rr|=0.
 \eess
 \end{theo}

{\bf Proof}. According to $0<c_0<\min\{\ud{c}_1,\,\ud{c}_2\}$ and (\ref{a.12}), there exist $0<\sigma_0\ll 1$ and $ t_\sigma\gg 1$ such that
 \bess
 c_\sigma:=c_0+\sigma<\min\{\ud{c}_1,\,\ud{c}_2\}, \ \ \forall\ 0<\sigma\le\sigma_0; \ \ s_1(t),\ s_2(t)>c_\sigma t, \ \ \forall\ t\ge t_\sigma.
 \eess

{\it Step 1}: Similar to the above, the estimate (\ref{a.4}) holds for $v$. For any given $0<\ep\ll1$, there exists $ t_1>0$ such that $v(t,x)<1+\ep$ in $[t_1,\yy)\times[0,\yy)$. Enlarging $ t_1$ if necessary, we may think $s_1(t_1)>\frac\pi 2(\frac{d_1}{r_1a_\ep})^{1/2}$, where $a_\ep=1-k(1+\ep)$.

\vskip 3pt {\it Step 2}: Let $(w_1,g_1)$ be the unique solution of (\ref{a.5}) with  $\Lambda=(t_1,d_1,r_1,a_\ep,\mu_1,s_1(t_1))$ and $w_0(x)=u(t_1,x)$.
Then $s_1(t)\ge g_1(t)$, $u(t,x)\ge w_1(t,x)$ in $D_{t_1}^{g_1}$ by the comparison principle. And, $g_1(\yy)=\yy$ since $g_1(t_1)=s_1(t_1)>\frac\pi 2(\frac{d_1}{r_1a_\ep})^{1/2}$. By use of \cite[Theorem 3.1]{ZhaoW15}, we get
 \bess
 \lim_{t\to\infty}(g_1(t)-c_\ep t)=\rho\in\mathbb{R}, \ \ \lim_{t\to\infty}\|w_1(t,x)- q_\ep(c_\ep t+\rho-x)\|_{L^\infty([0,g_1(t)])}=0,\eess
where $(q_\ep, c_\ep)$ is the unique solution of (\ref{a.2}) with $\gamma=(\mu_1,r_1a_\ep, r_1,d_1)$, i.e., $c_\ep=c_1(\mu_1,r_1a_\ep)$. Note $0<c_\sigma<\ud{c}_1$, we have $c_\ep>c_\sigma$ as $0<\ep\ll 1$.
Thus, $g_1(t)-c_\sigma t\to\yy$ and $\min_{[0,\,c_\sigma t]}(c_\ep t+\rho-x)\to\yy$ as $t\to\yy$.
Similar to the proof of (\ref{a.9}) we can derive
 \[\liminf_{t\to\yy}\min_{[0,\,c_\sigma t]}u(t,\cdot)\ge 1-k. \]
There exists $ t_2\gg 1$ such that
 \[s_2(t)>c_\sigma t, \ \ u(t,x)\ge 1-k-\ep:=b_\ep,\ \ \ \forall \ t\ge t_2, \ 0\le x\le c_\sigma t.\]

{\it Step 3}: Similar to the arguments of \cite[Lemma 2.1]{WZjdde}, $\limsup_{t\to\yy}v(t,x)\le 1-h b_\ep$ uniformly in $x\in[0,1]$. In particular,   $v(t,0)\le 1-h b_\ep+\ep$ in $[t_3,\yy)$ for some $ t_3> t_2$. Thus, $v$ satisfies
   \bess
 \left\{\begin{array}{lll}
  v_t-d_2v_{xx}\le r_2v(1-h b_\ep-v),\ \ &t\ge t_3, \ \ 0<x<c_\sigma t,\\[0.5mm]
 v(t,0)\le 1-h b_\ep+\ep, \ \ v<1+\ep, \ &t\ge t_3, \ 0<x\le c_\sigma t.
 \end{array}\right.
 \eess
We will show that $\limsup_{t\to\yy}\max_{[0,\,c_{\rho\sigma/2}t]}v(t,\cdot)\le 1-h b_\ep+\ep$, which leads to
 \bes
 \limsup_{t\to\yy}\max_{[0,\,c_{\sigma/2}t]}v(t,\cdot)\le 1-h(1-k):=\bar v_2\lbl{a.13}\ees
since $\ep>0$ is arbitrary. To do this, we choose  $0<\delta\ll 1$ and define
 \[\vp(t,x)=1-h b_\ep+\ep+h b_\ep{\rm e}^{\delta c_\sigma t_3}{\rm e}^{\delta(x-c_\sigma t)}, \ \ t\geq t_3, \ 0\le x\le c_\sigma t.\]
Evidently,
 \[\max_{[0,\,c_{\sigma/2} t]}\vp(t,\cdot)\le 1-h b_\ep+\ep+h b_\ep{\rm e}^{\delta c_\sigma t_3}{\rm e}^{-\delta\sigma t/2}\to 1-h b_\ep+\ep\]
as $t\to\yy$, and
 \[\vp(t,0)>1-h b_\ep+\ep, \ \ \vp(t,c_\sigma t)\ge 1+\ep, \ \ t\ge t_3; \ \
 \vp(t_3,x)\ge 1+\ep, \ \ 0\le x\le c_\sigma t_3.\]
Due to $0<\ep\ll1$, we can think of $1-h b_\ep+\ep>\frac 12[1-h(1-k)]$. 
It is easy to verify that 
  \[\vp_t-d_2\vp_{xx}\ge r_2\vp(1-h b_\ep-\vp), \ \ t\geq t_3, \ 0\le x\le c_\sigma t\]
provided that $\delta$ satisfies $\delta(c_\sigma+d_2\delta)\le \frac{r_2}2[1-h(1-k)]$. 
The comparison principle gives $v(t,x)\le\vp(t,x)$ for all $t\ge t_3$ and $0\le x\le c_\sigma t$. So, (\ref{a.13}) holds. We write $c_{\sigma/2}$ as $c_\sigma$ for the sake of writing. Then there exists $ t_4> t_3$ such that
 \[s_1(t)>c_\sigma t, \ \ v(t,x)\le \bar v_2+\ep:=\bar v_2^\ep<1, \ \ \ \forall \ t\ge t_4, \ 0\le x\le c_\sigma t.\]

{\it Step 4}: Because of $\bar v_2^\ep<1$, we have $c_1(\mu_1,r_1(1-k\bar v_2^\ep))>c_\sigma$. Take $0<\mu_1^*<\mu_1$ so that $c_1(\mu_1^*,r_1(1-k\bar v_2^\ep))=c_\sigma$. Let $(q_\sigma, c_\sigma)$ be the unique solution of (\ref{a.2}) with $\gamma=(\mu_1^*,r_1(1-k\bar v_2^\ep), r_1,d_1)$. Owing to $s_1(t)>c_\sigma t$ for all $t\ge t_4$, we can find a function $\tilde u\in C^2([0,c_\sigma t_4])$ satisfying $\tilde u(x)\le \min\{q_\sigma(c_\sigma t_4-x),\,u(t_4,x)\}$ in $[0,c_\sigma t_4]$ and
   \[\tilde u'(0)=\tilde u(c_\sigma t_4)=0, \ \tilde u(x)>0\ \ {\rm in}\, \
  [0,c_\sigma t_4).\]
Let $(w_2,g_2)$ be the unique solution of (\ref{a.5}) with $\Lambda=(t_4,d_1,r_1,1-k\bar v_2^\ep,\mu_1^*,c_\sigma t_4)$ and $w_0(x)=\tilde u(x)$.
Then, by use of \cite[Theorem 3.1]{ZhaoW15},
 \bess
 g_2(t)-c_\sigma t\to\rho\in\mathbb{R}, \ \ \|w_2(t,x)- q_\sigma(c_\sigma t+\rho-x)\|_{L^\infty([0,g_2(t)])}\to 0\eess
as $t\to\yy$. Define $z(t,x)=q_\sigma(c_\sigma t-x)$, $\eta(t)=c_\sigma t$. It is easy to verify that
 \bess\left\{\begin{array}{ll}
 z_t-d_1z_{xx}=r_1z(1-k\bar v_2^\ep-z), \ \ &t\ge  t_4, \ \ 0\le x\le\eta(t),\\[.2mm]
 -z_x(t,0)>0, \ \ z=0, \ \ \eta'(t)=-\mu_1^*z_x, \ \ &t\ge t_4, \ x=\eta(t),\\[.2mm]
 \eta(t_4)=g_2(t_4), \ \ z(t_4,x)\ge \tilde u(x)=w_2(t_4,x), \ \
   &0\le x\le c_\sigma t_4.
   \end{array}\right.\eess
By the comparison principle,
 \[g_2(t)\le\eta(t)=c_\sigma t, \ \ w_2(t,x)\le z(t,x)=q_\sigma(c_\sigma t-x) \ \ {\rm in} \ D_{t_4}^{g_2}.\]
Note that $g_2(t)\le c_\sigma t<s_1(t)$, $w_2(t,g_2(t))=0<u(t,g_2(t))$ in $[t_4,\yy)$, $w_2(t_4,x)\le u(t_4,x)$ in $[0,c_\sigma t_4]$,  and
 \[u_t-d_1u_{xx}\ge r_1u(1-k\bar v_2^\ep-u) \ \ {\rm in} \ D_{t_4}^{g_2}.\]
We have $u\ge w_2$ in $D_{t_4}^{g_2}$ by the comparison principle.
Similarly to Step 2, it can be derived that
 \[\liminf_{t\to\yy}\min_{[0,\,c_{\sigma/2}t]}u(t,\cdot)\ge 1-k\bar v_2= 1-k[1-h(1-k)]:=\ud{u}_2.\]

{\it Step 5}: Define
 \[\bar v_1=1,\ \ \ud{u}_1=1-k, \ \ \bar v_n=1-h\ud{u}_{n-1}, \ \ \ud{u}_n=1-k\bar v_n, \ \ n\ge 2.\]
Then $\ud{u}_n\to \frac{1-k}{1-hk}$, $\bar v_n\to\frac{1-h}{1-hk}$ as $n\to\yy$. Repeating the above process we can prove that
 \[\liminf_{t\to\yy}\min_{[0,\,c_0 t]}u(t,\cdot)\ge\ud{u}_n, \ \ \limsup_{t\to\yy}\max_{[0,\,c_0 t]}v(t,\cdot)\le\bar v_n, \ \ \forall\ n\geq 1.\]
Consequently,
 \[\liminf_{t\to\yy}\min_{[0,\,c_0 t]}u(t,\cdot)\ge(1-k)/(1-hk), \ \ \limsup_{t\to\yy}\max_{[0,\,c_0 t]}v(t,\cdot)\le(1-h)/(1-hk).\]

Similarly, we can show that
 \[\limsup_{t\to\yy}\max_{[0,\,c_0 t]}u(t,\cdot)\le(1-k)/(1-hk), \ \ \liminf_{t\to\yy}\min_{[0,\,c_0 t]}v(t,\cdot)\ge(1-h)/(1-hk).\]
The proof is complete. \ \ \ \fbox{}

Theorem \ref{th1.7} is an improvement of \cite[Theorem 6]{Wu15} in there the condition $hk<1/2$ is required.

\end{document}